\documentclass[a4paper,10pt]{article}
\usepackage{amssymb,amsbsy,amsmath,amsfonts,amscd}
\usepackage[T1]{fontenc}
\usepackage{textcomp}
\usepackage{graphicx}
\usepackage{color}
\usepackage{pdfsync}

\newtheorem{theorem}{Theorem}[section]
\newtheorem{lemma}[theorem]{Lemma}

\newtheorem{corollary}[theorem]{Corollary}

\begin{document}
\title{On problems in the calculus of variations \\
in increasingly elongated domains}

\author{Herv\'e Le Dret\\ Sorbonne Universit\'es,  UPMC Univ Paris 06, CNRS,\\Laboratoire Jacques-Louis Lions, Bo\^\i te courrier 187,\\
 75252 Paris Cedex 05, France.  Email: herve.le\_dret@upmc.fr\and
Amira Mokrane\\ Laboratoire d'\'equations aux d\'eriv\'ees partielles\\ non lin\'eaires et histoire des math\'ematiques, ENS, B.P. 92,\\ Vieux Kouba, 16050 Alger, Alg\'erie\\ and USTHB, Facult\'e des math\'ematiques, D\'epartement d'analyse,\\ Laboratoire d'analyse math\'ematique et num\'erique\\ des \'equations aux d\'eriv\'ees partielles, Bab Ezzouar, Alger, Alg\'erie.\\
 Email: mokrane\_amira3@yahoo.fr}

\maketitle

\begin{abstract}\noindent We consider minimization problems in the calculus of variations set
in a sequence of domains the size of which tends to infinity in certain directions and such that the data only depend on the coordinates in the directions that remain constant. We study the asymptotic behavior of minimizers in various situations and show that they converge in an appropriate sense toward minimizers of a related energy functional in the constant directions.
\end{abstract}

\section{Introduction}

In this article, we revisit the ``$\ell\to+\infty$'' problem in the context of the calculus of variations. This class of problems was introduced by Chipot and Rougirel in 2000, \cite{[CR1]}, see also the monograph by Chipot \cite{[C2]}, and has since then given rise to many works by several authors dealing with various elliptic and parabolic problems up to until recently. 

A prototypical $\ell\to+\infty$ problem is the following. Let  $\omega={]-}1,1[$, $\ell>0$ be a real number and $\Omega_{\ell}\subset{\mathbb R}^2$ the
rectangle ${]-}\ell,\ell[\times\omega$. We denote by $x_1$ the first variable in ${]-}\ell,\ell[$ and $x_2$ the second variable in $\omega$. Any function $f\in L^2(\omega)$ in the second variable gives rise to a function in two variables still denoted $f$ by setting $f(x_1,x_2)=f(x_2)$. We thus consider the two boundary value problems:
find $u_{\ell}$, a function in $(x_1,x_2)$, such that
$$\begin{cases}
           -\Delta u_{\ell}=f & \mbox{in } \Omega_{\ell},\\
           u_{\ell}=0 & \mbox{on } \partial\Omega_{\ell},
         \end{cases}
$$
and find $u_{\infty}$, a function in $x_2$, such that
$$\begin{cases}
           -\frac{d^2u_{\infty}}{dx_2^2}=f & \mbox{in } \omega,\\
           u_{\infty}=0 & \mbox{on }\partial\omega= \{-1,1\}.
         \end{cases}
$$
Now the function $u_\infty$ can also be considered as a function in two variables that is independent of $x_1$. In this case, it can be shown that, for any $\ell_0>0$, one has
$$u_{\ell}\rightarrow u_{\infty} \mbox{ in } H^1(\Omega_{\ell_0}) \hbox{ when }\ell\to+\infty,$$
hence the name of the problem. In other words, when the data does not depend on the elongated dimension, the solution of the above boundary value problem converges in some sense at finite distance to the solution of the corresponding boundary value problem posed in the non-elongated dimension when the elongation tends to infinity. 

The majority of works on the $\ell\to+\infty$ problem makes use of the boundary value problem itself, \emph{i.e.}, the PDE plus boundary condition. One exception to this rule are the recent papers {\cite{Chipot-Savitska,Chipot Mosjic Roy}}, in which the authors consider instead a sequence of problems in the calculus of variations posed on elongated domains, {see also \cite{[C2new]}. This is the approach we adopt here as well. 

Our main motivation for this is that certain models, such as nonlinear hyperelasticity, are naturally posed as problems in the calculus of variations for which no Euler-Lagrange equation, \emph{i.e.}, non underlying PDE even in a weak form, is available, see \cite{Ball}. Moreover, questions surrounding the Saint Venant principle in elasticity, see \cite{Mielke,Toupin}, are typically set in elongated domains, albeit in one direction only. Consequently, it makes sense to attempt dealing with some $\ell\to+\infty$ problems by using only energy minimization properties and no Euler-Lagrange equation whatsoever. 

We are however quite far from achieving the goal of treating nonlinear elasticity, since the approach that we develop below relies a lot on convexity, whereas convexity is not an appropriate hypothesis for nonlinear elasticity. We are nonetheless able to encompass a wide range of nonlinear energies, including the $p$-Lapla\-cian with some technical restrictions on the number of elongated dimensions with respect to the exponent $p$. Our hypotheses are weaker and our results are sometimes stronger than those of  \cite{Chipot Mosjic Roy}. The techniques are somewhat different too, with an emphasis here on weak convergence and weak lower semicontinuity techniques, and reliance on such classical techniques as the De Giorgi slicing method which are not dependent on convexity. As a general rule, we try to make as little use of convexity as we can at any given point. 

Let us describe our results a little more precisely. We consider bounded open subsets $\Omega_\ell$ of ${\mathbb R}^n$  which are Cartesian products of the form $\ell\omega'\times\omega''$, with $\omega'\subset {\mathbb R}^r$ and $\omega''\subset{\mathbb R}^{n-r}$, with $1\le r\le n-1$. We let $x=(x',x'')$ with $x'\in {\mathbb R}^r$ being the elongated variable and $x''\in {\mathbb R}^{n-r}$ the non-elongated variable. Likewise, for a scalar-valued function $v\colon\Omega_\ell\to{\mathbb R}$, we decompose the gradient $\nabla v=(\nabla'v,\nabla''v)$ with obvious notation.
 
We consider an energy density $F\colon{\mathbb R}^n\to{\mathbb R}$ and a function $f$ on $\omega''$, and introduce the minimization problem of finding $u_\ell\in W^{1,p}_0(\Omega_\ell)$ such that $J_\ell(u_\ell)=\inf_{v\in W^{1,p}_0(\Omega_\ell)}J_\ell(v)$ where
$$
J_\ell(v)=\int_{\Omega_\ell}\bigl(F(\nabla v(x))-f(x'')v(x)\bigr)\,dx.$$
We assume that $F$ has $p$-growth, $p$-coerciveness and is convex. In particular, there is no assumption of strict convexity or uniform strict convexity made on $F$. 

We then introduce $F''\colon {\mathbb R}^{n-r}\to{\mathbb R}$ by letting $F''(\xi'')=F(0,\xi'')$, again with obvious notation. Of course, $F''$ is convex, has $p$-growth and $p$-coerciveness and the minimization problem of finding $u_\infty\in W^{1,p}_0(\omega'')$ such that $J_\infty(u_\infty)=\inf_{v\in W^{1,p}_0(\omega'')}J_\infty(v)$ where
$$
J_\infty(v)=\int_{\omega''}\bigl(F''(\nabla'' v(x''))-f(x'')v(x'')\bigr)\,dx'',$$
admits solutions. It turns out that, under additional hypotheses, this problem is the ``$\ell\to+\infty$'' limit of the family of minimization problems under consideration.

These hypotheses include appropriate growth and coerciveness hypotheses on the function $G\colon {\mathbb R}^n\to{\mathbb R}$, $G(\xi)=F(\xi)-F''(\xi'')$, of the form
$$
\forall \xi\in{\mathbb R}^{n},
\lambda(|\xi'|^p+k|\xi''|^{p-k}|\xi'|^k)\leq G(\xi)\leq\Lambda(|\xi'|^{p}+k|\xi''|^{p-k}|\xi'|^k),
$$
for some $0<\lambda\le\Lambda$ and $0\le k<p$. Depending on the case, there is no more additional hypothesis (for $k=0$), or a hypothesis of strict convexity of $F''$, or a hypothesis of uniform strict convexity of $F''$ (for $k>0)$.

The results are a ``$\ell\to+\infty$'' convergence in the weak sense for $k=0$ when $r<p$, sharpened to strong sense when $F''$ is furthermore assumed to be strictly convex, and a strong ``$\ell\to+\infty$'' convergence for $k>0$ when $r\leq kp/(p-k)$. In the case of the $p$-Laplacian, $p>2$, we thus obtain strong ``$\ell\to+\infty$'' convergence when $r<2p/(p-2)$, {see also \cite{Xie}}.

In addition, in the case $k=0$, if we assume that $F''$ is uniformly strictly convex, we obtain strong convergence  at an exponential rate without any restriction on $r$.  This includes the known behavior of the $2$-Laplacian in the ``$\ell\to+\infty$'' context.

We conclude the article with a few comments and perspectives on the vectorial case, in connection with  nonlinear elasticity in particular.

\section{Statement of the problem}

We consider two bounded open sets $\omega'\subset{\mathbb R}^r$ with $0\in\omega'$ and $\omega'$ is starshaped with respect to $0$, and $\omega''\subset{\mathbb R}^{n-r}$ with $n>r\ge 1$. Let $\ell >0$ and set
\begin{equation}\label{omgl}
\omega'_\ell=\ell\omega'\text{ and }\Omega_\ell=\omega'_\ell\times\omega''\subset{\mathbb R}^n.
\end{equation}

Points $x$ in $\Omega_\ell$ will be denoted by $x=(x',x'')$ with $x'=(x_1,x_2,\ldots,x_r)\in\omega'_\ell$ and $x''=(x_{r+1},\ldots,x_n)\in \omega''$. Likewise, vectors $\xi$ in ${\mathbb R}^n$ will be decomposed as $\xi=(\xi',\xi'')$, with $\xi'\in {\mathbb R}^r$ and $\xi''\in {\mathbb R}^{n-r}$.

Note that because of the starshaped assumption, we have $\Omega_\ell\subset\Omega_{\ell'}$ as soon as $\ell\le\ell'$ and we are thus dealing with a ``growing' family of open sets.
We make an additional regularity hypothesis on $\omega'$, which is as follows. Define first the gauge function of $\omega'$ as
$$g(x')=\inf\{t\in{\mathbb R}_+^*; x'/t\in \omega'\}.$$
Since $\omega'$ is starshaped and bounded, this is well defined, $\omega'_\ell=\{x'; g(x')<\ell\}$, and there exists $0<R_1<R_{2}$  such that $R_1|x'|\le g(x')\le R_2|x'|$ for all $x'\in{\mathbb R}^r$.

Now we assume that $\omega'$ is such that $g$ is a Lipschitz function with Lipschitz constant $K$. By Rademacher's theorem, this implies that $g$ is almost everywhere differentiable, with $|\nabla'g(x')|\le K$ a.e. Moreover, it is known that $g$ then belongs to $W^{1,\infty}_{\rm loc}({\mathbb R}^r)$ and that its almost everywhere derivatives equal its distributional derivatives. This is true for example if $\omega'$ is convex. This regularity hypothesis is for convenience only: we use $g$ to build cut-off functions inside the domains, and not up to the  boundary. It should be quite clear that our results can be rewritten in order to accommodate arbitrary open sets $\omega'$.

We are interested in a sequence of problems
in the calculus of variations ${\cal P}_\ell$ of the form
\begin{equation}\label{pl}
J_\ell(u_\ell)=\inf_{v\in W^{1,p}_0(\Omega_\ell)}J_\ell(v),
\end{equation}
with $u_\ell\in W^{1,p}_0(\Omega_\ell)$ and
\begin{equation}\label{funjl}
J_\ell(v)=\int_{\Omega_\ell}\bigl[F(\nabla v(x))-f''(x'')v(x)\bigr]\,dx.
\end{equation}
where $f''\in L^{p'}(\omega'')$, $\frac1p+\frac1{p'}=1$, is a given function.
Observe that the term corresponding to the force term for this problem only depends on the ``non-elongated'' variable $x''$ so that it is reasonable to expect that $u_\ell$ behaves
 as a function mostly in $x''$ in the limit $\ell\to+\infty$, in a sense made precise below. We could also consider more general semilinear force terms of the form $h(x'',v)$  satisfying appropriate growth and convexity assumptions, but we stick here with a linear term for simplicity.

We assume that the energy density $F\colon{\mathbb R}^n\to{\mathbb R}$ is convex. We let
\begin{equation}\label{structure de F}
\begin{array}{rcl}
F''\colon {\mathbb R}^{n-r}&\to&{\mathbb R}\\
\xi''&\mapsto&F(0,\xi'')
\end{array}
\quad\hbox{and}\quad
\begin{array}{rcl}
G\colon {\mathbb R}^n&\to&{\mathbb R}\\
\xi&\mapsto&F(\xi)-F''(\xi'')
\end{array}
\end{equation}
so that
\begin{equation}\label{f1f2}
F(\xi',\xi'')=F''(\xi'')+G(\xi',\xi''),
\end{equation}
and
$F''$ is convex. These functions are assumed to satisfy the following coerciveness and growth hypotheses
\begin{align}
&\forall \xi\in{\mathbb R}^{n},
\lambda|\xi|^p\leq F(\xi)\leq\Lambda(|\xi|^p+1),\label{growth1}\\
&\forall \xi\in{\mathbb R}^{n},
\lambda(|\xi'|^p+k|\xi''|^{p-k}|\xi'|^k)\leq G(\xi)\leq\Lambda(|\xi'|^{p}+k|\xi''|^{p-k}|\xi'|^k),\label{growth3}
\end{align}
for some $0<\lambda\le\Lambda$, $p>1$ and $0\le k< p$.\footnote{Note that  $k=p$ yields the same hypothesis as $k=0$.} Here, for $\xi\in{\mathbb R}^d$, $|\xi|$ denotes the canonical Euclidean norm of $\xi$ in ${\mathbb R}^d$.

Clearly, condition \eqref{growth1} implies the similar condition
\begin{equation}
\forall \xi''\in{\mathbb R}^{n-r},
\lambda|\xi''|^p\leq F''(\xi'')\leq\Lambda(|\xi''|^p+1),\label{growth2}
\end{equation}
for $F''$.

Energy densities of the form above include that associated with the $p$-Laplacian for $p\ge 2$. Indeed, in this case, $F(\xi)=\frac1p|\xi|^p=\frac1p(|\xi'|^2+|\xi''|^2)^{p/2}$ and we can take $k=2$ for $p>2$, or $k=0$ for $p=2$. Another simple energy density that is covered by our analysis is $F(\xi)=\frac1p (|\xi'|^p+|\xi''|^p)$ or more generally energies of the form $F(\xi)=F'(\xi')+F''(\xi'')$, with appropriate hypotheses on $F'$ and $F''$. Here, assuming without loss of generality that $F'(0)=0$, we have $G(\xi',\xi'')=F'(\xi')$ and we can take $k=0$.

In addition to the above growth and coerciveness hypotheses, which obviously imply  that problem ${\cal P}_\ell$ has at least one solution $u_\ell$, we assume that $F''$ is uniformly strictly convex for $k>0$, in the sense that there exists a constant $\beta>0$ such that for all $\xi''$, $\zeta''\in {\mathbb R}^{n-r}$ and all $\theta,\mu\in [0,1]$ with $\theta+\mu=1$, we have
\begin{equation}\label{uniformite stricte}
F''(\theta\xi''+\mu\zeta'')\le \theta F''(\xi'')+\mu F''(\zeta'')-k\beta\theta\mu(\theta^{p-1}+\mu^{p-1})|\xi''-\zeta''|^p.
\end{equation}
see for instance \cite{[A],Evans,[J-N]}. The $p$-Laplacian for $p> 2$, $k=2$, satisfies this hypothesis (the $2$-Laplacian satisfies the alternate hypothesis \eqref{uniformite stricte k=0} that will be used later on in Section 5).
Note that when $k=0$, the hypothesis becomes redundant, and there is actually no requirement of even strict convexity, let alone uniform strict convexity, of $F''$ in this case.

We now introduce our candidate limit problem
${\cal P}_\infty$ as that of finding $u_\infty\in W^{1,p}_0(\omega'')$ such that
\begin{equation}\label{pinfty}
J_\infty(u_\infty)=\inf_{v\in W^{1,p}_0(\omega'')}J_\infty(v),
\end{equation}
with
\begin{equation}\label{funjinfty}{}
J_{\infty}(v)=\int_{\omega''}\bigl[F''(\nabla''v(x''))-f''(x'')v(x'')\bigr]\,dx''.
\end{equation}
It also clear that problem ${\cal P}_\infty$ has at least one solution $u_\infty$.

Here and in the sequel, we use the following notational device
$$\nabla'=(\partial_1,\ldots,\partial_r),\quad\nabla''=(\partial_{r+1},\ldots,\partial_n),$$
that we apply indifferently to functions defined either on $\Omega_\ell$ or  on $\omega''$. For brevity, we refer to $\nabla'$ as the ``horizontal'' part of the gradient and to $\nabla''$ as the ``vertical'' part of the gradient.

We want to study the asymptotic behavior of $u_\ell$ when $\ell\to+\infty$ and compare it with a minimizer $u_\infty$ of the $n-r$ dimensional vertical problem ${\cal P}_\infty$. Actually, our goal is to show that the former converges to the latter in a sense that will be explained later on.

\section{Preliminary estimates}

We first give several estimates that we will use in the proofs of our convergence
results. The first estimate follows immediately from Poincar\'e's inequality.
\begin{lemma}\label{lem1}
There exists a constant $c_1=c_1(\omega'')$ independent of $\ell$ such that for all
 $v\in W^{1,p}(\Omega_{\ell})$ whose trace vanish on
$\omega'_\ell\times\partial\omega''$, we have
\begin{equation}\label{poincare}
\|v\|_{L^p(\Omega_{\ell})}\leq c_1 \|\nabla''
v\|_{L^p(\Omega_{\ell};{\mathbb R}^{n-r})}.
\end{equation}
\end{lemma}

Let us now give a first, coarse estimate of $u_\ell$.

\begin{lemma}\label{lem2}
There exists a constant $c_2$ independent of $\ell$, such that
\begin{equation}\label{estmtul}
\int_{\Omega_{\ell}}|\nabla u_{\ell}|^p\,dx\le c_2 \ell^r.
\end{equation}
\end{lemma}

\noindent{\bf Proof.}\enspace{}
Let us take $v=0$ as a test-function in problem \eqref{pl}. It follows that
$$\int_{\Omega_{\ell}} F(\nabla u_{\ell}(x))\,dx\le\int_{\Omega_{\ell}}
f''(x'')u_{\ell}(x)\,dx+A\ell^r.$$
where $A=F(0){\cal L}^r(\omega'){\cal L}^{n-r}(\omega'')$ does not depend on $\ell$ (${\cal L}^d$ denotes the $d$-dimensional Lebesgue measure). By H\"older's inequality and the coerciveness assumption
\eqref{growth1}, it follows that
\begin{align*}
\int_{\Omega_{\ell}} |\nabla u_{\ell}(x)|^p\,dx &\le
\frac{1}{\lambda}\Bigl(\int_{\Omega_{\ell}}|f''(x'')|^{p'}\,dx\Bigr)^{1/p'}
\Bigl(\int_{\Omega_{\ell}}|u_{\ell}(x)|^p\,dx \Bigr)^{1/p} +\frac {A}{\lambda}\ell^r\\
 &\le
 \frac{B}{\lambda}
\ell^{r/p'}
\|\nabla''
u_{\ell}\|_{L^p(\Omega_{\ell};{\mathbb R}^{n-r})}+\frac {A}{\lambda}\ell^r,
\end{align*}
with $B=c_1\|f''\|_{L^{p'}(\omega'')}{\cal L}^r(\omega')^{1/p'}$, which does not depend on $\ell$. Consequently, we obtain an estimate of the form
\begin{equation}\label{estimation grossiere}
\|\nabla
u_{\ell}\|^p_{L^p(\Omega_{\ell};{\mathbb R}^{n})} \le
 C
\ell^{r/p'}
\|\nabla
u_{\ell}\|_{L^p(\Omega_{\ell};{\mathbb R}^{n})}+D\ell^r,
\end{equation}
where $C$ and $D$ are constants that do not depend on $\ell$. Let us set $X=\ell^{-r/p}\|\nabla
u_{\ell}\|_{L^p}$. Estimate \eqref{estimation grossiere} now reads
$$X^p\le CX+D,$$
so that there exists $c_2$ depending only on $C$ and $D$ such that $X\le c_2^{1/p}$, which completes the proof.\hfill$\square$\par\medbreak

We now recall an elementary estimate similar to what can be found
 in \cite{[G2]} for instance.
\begin{lemma}\label{lem3}
Let $h(t)$ a nonnegative bounded function defined on an interval
$[\tau_0,\tau_1]$, $\tau_0\geq 0$. Suppose that for $\tau_0\leq t<
s\leq \tau_1$, we have
$$h(t)\leq \theta h(s)+C(s-t)^{-\nu_1}+D(s-t)^{-\nu_2},$$
where $C,D,\nu_1,\nu_2,\theta$ are nonnegative constants with
$0\leq \theta< 1$. Then, for all $\tau_0\le t<s\le\tau_1$, we have
$$h(t)\le
c(C(s-t)^{-\nu_1}+D(s-t)^{-\nu_2}),$$
where $c$ is a constant that only depends on
$\nu_1$, $\nu_2$ and $\theta$.
\end{lemma}

\noindent{\bf Proof.}\enspace{} If we have two sequences of nonnegative numbers $a_i$ and $b_i$ such that $a_i\le \theta a_{i+1}+b_{i+1}$, it follows by induction that $a_0\le \theta^i a_i+\sum_{j=0}^{i-1}\theta^jb_{j+1}$. We apply this remark to the sequences $a_i=h(t_i)$ and $b_{i+1}=C(t_{i+1}-t_i)^{-\nu_1}+D(t_{i+1}-t_i)^{-\nu_2}$, where $t_i=t+(1-\sigma^i)(s-t)$, $0<\sigma<1$ to be chosen later on, is an increasing sequence in $[\tau_0,\tau_1]$ such that $t_0=t$. This yields the estimate
$$h(t)\le \theta^ih(t_i)+\frac{C}{(s-t)^{\nu_1}}(1-\sigma)^{-\nu_1}\sum_{j=0}^{i-1}\biggl(\frac\theta{\sigma^{\nu_1}}\biggr)^j+\frac{D}{(s-t)^{\nu_2}}(1-\sigma)^{-\nu_2}\sum_{j=0}^{i-1}\biggl(\frac\theta{\sigma^{\nu_2}}\biggr)^j .$$
We now choose $\sigma<1$ in such a way that $\frac\theta{\sigma^{\nu_1}}<1$ and $\frac\theta{\sigma^{\nu_2}}<1$, and conclude by letting $i\to+\infty$, remembering that $h(t_i)$ is bounded.\hfill$\square$\par\medbreak

Next, we estimate the horizontal part of the gradient of $u_\ell$ in $L^p(\Omega_{\ell_0})$ in terms of $\ell$, $\ell_0$, $u_\ell$ and a minimizer $u_\infty$ of the vertical problem ${\cal P}_\infty$.
\begin{theorem}\label{estimation de base} There exists a constant $c_3$ independent of all the other quantities such that, for all $0<t<s\le\ell$ and all minimizers $u_\infty$ of the vertical problem, we have
\begin{multline}\label{l'estimation en question usc}
\|\nabla'u_\ell\|^p_{L^p(\Omega_{t};{\mathbb R}^r)}+k\|\nabla''(u_\ell-u_\infty)\|^p_{L^p(\Omega_t;{\mathbb R}^{n-r})}\\\le \frac{\delta c_3}{(s-t)^p}\|\nabla''(u_\ell-u_\infty)\|^p_{L^p(\Omega_s\setminus\Omega_t;{\mathbb R}^{n-r})}\\+
\frac{c_3k}{(s-t)^{kp/(p-k)}}\left\{\|\nabla''u_{\ell}\|^p_{L^p(\Omega_s\setminus\Omega_t;{\mathbb R}^{n-r})}+(1-\delta)\|\nabla''(u_\ell-u_\infty)\|^p_{L^p(\Omega_s\setminus\Omega_t;{\mathbb R}^{n-r})}\right\}.
\end{multline}
where $\delta=1$ if $0\le k\le p/2$, $\delta=0$ otherwise.
\end{theorem}

\noindent{\bf Proof.}\enspace{} We first define a family of cut-off functions as follows. For all $0<t<s\le \ell$, we set
$$\rho_{s,t}(x')=\frac1{s-t}\min\{(s-g(x'))_+,s-t\}.$$
By the definition of the gauge function, we see that $\rho_{s,t}\equiv 0$ on $\omega'_\ell\setminus\omega'_s$, $\rho_{s,t}\equiv 1$ on $\omega'_t$ and $0\le \rho_{s,t}\le 1$.
By our regularity assumption on $\omega'$, $\rho_{s,t}$ is Lipschitz and such that
$$\nabla'\rho_{s,t}(x')=-\frac1{s-t}\nabla'g(x'){\bf 1}_{\omega'_s\setminus\omega'_t}(x'),$$
so that we can estimate
\begin{equation}\label{gradient cutoff}
|\nabla'\rho_{s,t}(x')|\le\frac K{s-t}{\bf 1}_{\omega'_s\setminus\omega'_t}(x').
\end{equation}

We pick a number $0<\alpha<1$ and then set
\begin{equation}\label{test1}
v_1(x)=(1-\alpha\rho_{s,t}(x'))u_{\ell}(x)+\alpha\rho_{s,t}(x')u_{\infty}(x''),
\end{equation}
and
\begin{equation}\label{test2}
v_2(x)=(1-\alpha\rho_{s,t}(x'))u_{\infty}(x'')+\alpha\rho_{s,t}(x')u_{\ell}(x).
\end{equation}
Clearly,  $v_1$ belongs to $W^{1,p}_0(\Omega_\ell)$ and is thus a suitable test-function for problem ${\cal P}_\ell$, hence
\begin{equation}\label{test v1}
\int_{\Omega_\ell}\bigl[F(\nabla u_\ell(x))-f''(x'')u_\ell(x)\bigr]\,dx\le\int_{\Omega_\ell}\bigl[F(\nabla v_1(x))-f''(x'')v_1(x)\bigr]\,dx.
\end{equation}

Next we note that, owing to the embedding $W^{1,p}_0(\Omega_\ell)\hookrightarrow L^p(\omega'_\ell; W^{1,p}_0(\omega''))$, $v_2$ is suitable test-function for problem ${\cal P}_\infty$ for almost all $x'$, hence
\begin{multline}\label{test v2}
\int_{\omega''}\bigl[F''(\nabla''u_\infty(x''))-f''(x'')u_\infty(x'')\bigr]\,dx''\\
\le\int_{\omega''}\bigl[F''(\nabla''v_2(x',x''))-f''(x'')v_2(x',x'')\bigr]\,dx''.
\end{multline}
Integrating estimate \eqref{test v2} over $\omega'_\ell$, we obtain
\begin{multline}\label{test v2 bis}
\int_{\Omega_\ell}\bigl[F''(\nabla''u_\infty(x''))-f''(x'')u_\infty(x'')\bigr]\,dx\\
\le\int_{\Omega_\ell}\bigl[F''(\nabla''v_2(x))-f''(x'')v_2(x)\bigr]\,dx.
\end{multline}

We add estimates \eqref{test v1} and \eqref{test v2 bis} together and note that all the terms involving $f''$ cancel out since $v_1+v_2=u_\ell+u_\infty$. Therefore,
\begin{equation}\label{addition tests}
\int_{\Omega_\ell}\bigl[F(\nabla u_\ell(x))+F''(\nabla''u_\infty(x''))\bigr]\,dx
\le \int_{\Omega_\ell}\bigl[F(\nabla v_1(x))+F''(\nabla''v_2(x))\bigr]\,dx.
\end{equation}
We observe that $v_1=u_\ell$ and $v_2=u_\infty$ on $\Omega_\ell\setminus\Omega_s$, so that estimate \eqref{addition tests} boils down to
\begin{align}
\int_{\Omega_s}\bigl[F(\nabla u_\ell(x))+F''(\nabla''u_\infty(x''))\bigr]\,dx
&\le \int_{\Omega_s\setminus\Omega_t}\bigl[F(\nabla v_1(x))+F''(\nabla''v_2(x))\bigr]\,dx\nonumber\\
&\qquad+\int_{\Omega_t}\bigl[F(\nabla v_1(x))+F''(\nabla''v_2(x))\bigr]\,dx.\label{addition tests 2}
\end{align}
The left-hand side of \eqref{addition tests 2}
 can be rewritten as
\begin{multline}\label{addition tests gauche reecrite}
\int_{\Omega_s}\bigl[F(\nabla u_\ell(x))+F''(\nabla''u_\infty(x''))\bigr]\,dx\\=\int_{\Omega_s}\bigl[G(\nabla u_\ell(x))+F''(\nabla''u_\ell(x))+F''(\nabla''u_\infty(x''))\bigr]\,dx.
\end{multline}

Let $I_1$ and $I_2$ be the first and second integrals in the right-hand side of \eqref{addition tests 2}. To estimate $I_1$, we just use the convexity of $F''$, since the vertical gradients of $v_1$ and $v_2$ are convex combinations of the vertical gradients of $u_\ell$ and $u_\infty$,
\begin{align}
I_1&=\int_{\Omega_s\setminus\Omega_t}\bigl[G(\nabla v_1(x))+F''(\nabla'' v_1(x))+F''(\nabla''v_2(x))\bigr]\,dx\nonumber\\
&\le\int_{\Omega_s\setminus\Omega_t}\bigl[G(\nabla v_1(x))+F''(\nabla'' u_\ell(x))+F''(\nabla''u_\infty (x))\bigr]\,dx.\label{estimation I1}
\end{align}
To estimate $I_2$, we note that $v_1=(1-\alpha)u_\ell+\alpha u_\infty$ and $v_2=\alpha u_\ell+ (1-\alpha)u_\infty$ on $\Omega_t$, thus owing to the convexity of $F$ and the uniform convexity \eqref{uniformite stricte} of $F''$,
\begin{align}
I_2&\le \int_{\Omega_t}\bigl[(1-\alpha)F(\nabla u_\ell)+\alpha F(\nabla u_\infty)+(1-\alpha)F''(\nabla''u_\infty)+\alpha F''(\nabla''u_\ell)\nonumber\\
&\qquad\qquad\qquad\qquad\qquad\qquad\qquad\qquad\qquad\qquad -k\gamma|\nabla''(u_\infty-u_\ell)|^p\bigr]\,dx\nonumber\\
&= \int_{\Omega_t}\bigl[(1-\alpha)G(\nabla u_\ell)+F''(\nabla''u_\ell)+F''(\nabla'' u_\infty) -k\gamma|\nabla''(u_\infty-u_\ell)|^p\bigr]\,dx,\label{estimation I2}
\end{align}
for some $\gamma>0$. Putting estimates \eqref{addition tests 2}, \eqref{estimation I1}, \eqref{estimation I2} and equation \eqref{addition tests gauche reecrite} together, we obtain
$$\int_{\Omega_s\setminus\Omega_t}G(\nabla u_\ell)\,dx+\int_{\Omega_t}\bigl[\alpha G(\nabla u_\ell)+k\gamma|\nabla''(u_\infty-u_\ell)|^p\bigr]\,dx
\le\int_{\Omega_s\setminus\Omega_t}G(\nabla v_1(x))\,dx,$$
which, upon using the coerciveness hypothesis \eqref{growth3}, yields
\begin{multline}\label{premiere estimation serieuse}
a\int_{\Omega_t}\bigl[(|\nabla'u_\ell|^p+k|\nabla''u_\ell|^{p-k}|\nabla'u_\ell|^k)+k|\nabla''(u_\infty-u_\ell)|^p\bigr]\,dx
\\\le\int_{\Omega_s\setminus\Omega_t}G(\nabla v_1(x))\,dx,
\end{multline}
where $a>0$ is a small generic constant that only depends on the other constants involved.

We now focus on estimating the right-hand side of \eqref{premiere estimation serieuse}. We have
\begin{equation}\label{gradients tests}
\left\{\begin{aligned}
\nabla'v_1&=(1-\alpha\rho_{s,t})\nabla'u_{\ell}+\alpha\nabla'\rho_{s,t}(u_{\infty}-u_{\ell}),\\
\nabla''v_1&=(1-\alpha\rho_{s,t})\nabla''u_{\ell}+\alpha\rho_{s,t}\nabla''u_{\infty}.
\end{aligned}\right.
\end{equation}

Based on \eqref{gradients tests} and the definition of $\rho_{s,t}$, we have the following estimates for any exponent $q$:
\begin{equation}\label{gradients tests estimes}
\left\{\begin{aligned}
|\nabla'v_1|^q&\le2^{q-1}|\nabla'u_{\ell}|^q+2^{q-1}\frac{K^q}{(s-t)^q}|u_{\infty}-u_{\ell}|^q,\\
|\nabla''v_1|^{q}&\le2^{q-1}|\nabla''u_{\ell}|^{q}+2^{q-1}|\nabla''(u_{\infty}-u_{\ell})|^{q}.
\end{aligned}\right.
\end{equation}
We will use exponents $q=p$ and $q=k$ for the first line and $q=p-k$ for the second line. Due to the growth hypothesis \eqref{growth3}, we have
\begin{multline}\label{estimation de Gv1}
G(\nabla v_1)\le A\biggl(|\nabla'u_{\ell}|^p+\frac{1}{(s-t)^p}|u_{\infty}-u_{\ell}|^p\\+k\bigl(|\nabla''u_{\ell}|^{p-k}+|\nabla''(u_{\infty}-u_{\ell})|^{p-k}\bigr)\Bigl(|\nabla'u_{\ell}|^k+\frac{1}{(s-t)^k}|u_{\infty}-u_{\ell}|^k\Bigr)\biggr),
\end{multline}
where $A$ is a large generic constant that only depends on the other constants involved. For $k\ge 1$, three of the four product terms that appear need to be estimated. For this purpose, we will use Young's inequality in the following  form
$$a^kb^{p-k}\le \frac{k}{p}a^p+\frac{p-k}{p}b^{p}$$
for $a,b\ge 0$ (recall that $p> k$). We thus obtain
\begin{multline}\label{estimation de Gv1 2}
G(\nabla v_1)\le A\biggl(|\nabla'u_{\ell}|^p+\frac{1}{(s-t)^p}|u_{\infty}-u_{\ell}|^p\\
+k\Bigl(|\nabla''u_{\ell}|^{p-k}|\nabla'u_{\ell}|^k+|u_{\infty}-u_{\ell}|^p+\frac{1}{(s-t)^{kp/(p-k)}}|\nabla''u_{\ell}|^{p}
+|\nabla''(u_{\infty}-u_{\ell})|^{p}\Bigr)\biggr),
\end{multline}
where $A$ is another generic constant. We integrate this inequality over $\Omega_s\setminus\Omega_t$ and use Poincar\'e's inequality in the vertical variables to obtain
\begin{multline}\label{estimation du membre de droite}
\int_{\Omega_s\setminus\Omega_t}G(\nabla v_1)\,dx\le A\int_{\Omega_s\setminus\Omega_t}\bigl(|\nabla'u_{\ell}|^p+k|\nabla''u_{\ell}|^{p-k}|\nabla'u_{\ell}|^k+k|\nabla''(u_{\infty}-u_{\ell})|^{p}\bigr)\,dx\\
+\frac{A}{(s-t)^p}\|\nabla''(u_{\infty}-u_{\ell})\|^p_{L^p(\Omega_s\setminus\Omega_t)}
+\frac{Ak}{(s-t)^{kp/(p-k)}}\|\nabla''u_{\ell}\|^p_{L^p(\Omega_s\setminus\Omega_t)},
\end{multline}
with $A$ yet another generic constant.

We now consider two different cases. First, for $0\le k\le p/2$, 
let us set
$$h(t)=\int_{\Omega_{t}}\bigl(|\nabla' u_{\ell}|^p+k|\nabla'' u_{\ell}|^{p-k}|\nabla' u_{\ell}|^k+k|\nabla''(u_{\infty}-u_{\ell})|^{p}\bigr)\,dx.$$
Inequalities \eqref{premiere estimation serieuse} and \eqref{estimation du membre de droite} may be rewritten as
\begin{multline}\label{inegalite de g}
h(t)\le \theta h(s)
+\frac{1}{(s-t)^{p}}\|\nabla''(u_{\infty}-u_{\ell})\|^p_{L^p(\Omega_s\setminus\Omega_t)}\\
+\frac{k}{(s-t)^{kp/(p-k)}}\|\nabla''u_{\ell}\|^p_{L^p(\Omega_s\setminus\Omega_t)}
,
\end{multline}
with $\theta=\frac A{A+a}\in {]}0,1[$. Let $t\le t_1<s_1\le s$.
We  invoke Lemma \ref{lem3}, with $\nu_1=p$, $C=\|\nabla''(u_{\infty}-u_{\ell})\|^p_{L^p(\Omega_s\setminus\Omega_t)}$, $\nu_2=kp/(p-k)$, $D= k\|\nabla''u_{\ell}\|^p_{L^p(\Omega_s\setminus\Omega_t)}$, to conclude that
\begin{equation}\label{inegalite de g corrigee}
h(t_1)\le
c(C(s_1-t_1)^{-\nu_1}+D(s_1-t_1)^{-\nu_2}).
\end{equation}

The result follows in this case  by letting $t_1\to t$ and $s_1\to s$ since the constant $c$ only depends on $\nu_1$, $\nu_2$ and $\theta$, and $h$ is continuous (recall that $\delta=1$).

Now the second case is when $p/2<k<p$. Estimate \eqref{estimation du membre de droite} still holds true, but we now use Young's inequality once more in the form
$$a^wb^{p-w}\le \frac{w}{p}a^p+\frac{p-w}{p}b^{p}$$
with $w=\frac{p(2k-p)}k$ to deduce that
$$\frac1{(s-t)^p}=1^w\biggl(\frac1{(s-t)^{p/(p-w)}}\biggr)^{p-w}\le \frac{p-k}k+\frac{2k-p}k\frac1{(s-t)^{kp/(p-k)}},$$
so that we can actually write
\begin{equation}\label{inegalite de g bis}
h(t)\le \theta h(s)
+\frac{k}{(s-t)^{kp/(p-k)}}\bigl(\|\nabla''(u_{\infty}-u_{\ell})\|^p_{L^p(\Omega_s\setminus\Omega_t)}
+\|\nabla''u_{\ell}\|^p_{L^p(\Omega_s\setminus\Omega_t)}\bigr)
,
\end{equation}
with the same function $h$, but with another value for $\theta$, which we do not write here. We conclude as before with Lemma \ref{lem3} and the first constant $C=0$ for instance.\hfill$\square$\par\medbreak

The following is an immediate consequence of the previous estimate.
\begin{corollary}\label{Caccio global}
We have, for all $\ell\ge \ell_0$,
  \begin{multline}\label{Caccio usc}
\|\nabla'u_\ell\|^p_{L^p(\Omega_{\ell_0};{\mathbb R}^r)}+k\|\nabla''(u_\ell-u_\infty)\|^p_{L^p(\Omega_{\ell_0};{\mathbb R}^{n-r})}\\\le \frac{\delta c_3}{(\ell-\ell_0)^p}\|\nabla''(u_\ell-u_\infty)\|^p_{L^p(\Omega_\ell;{\mathbb R}^{n-r})}\\+\frac{c_3k}{(\ell-\ell_0)^{kp/(p-k)}}\left\{\|\nabla''u_{\ell}\|^p_{L^p(\Omega_\ell;{\mathbb R}^{n-r})}+(1-\delta)\|\nabla''(u_\ell-u_\infty)\|^p_{L^p(\Omega_\ell;{\mathbb R}^{n-r})}\right\},
\end{multline}
where $\delta=1$ if $0\le k\le p/2$, $\delta=0$ otherwise.
\end{corollary}

\noindent{\bf Proof.}\enspace{} Indeed, we take $s=\ell$, $t=\ell_0$ and notice that $\Omega_\ell\setminus\Omega_{\ell_0}\subset\Omega_\ell$.\hfill$\square$\par\medbreak

Let us remark that if $k=0$ and there is actually no strict convexity assumption made on $F''$, \emph{i.e.}, $F''$ may well be not strictly convex, the previous result boils down to 
$$
\|\nabla'u_\ell\|^p_{L^p(\Omega_{\ell_0};{\mathbb R}^r)}\le \frac{c_3}{(\ell-\ell_0)^p}\|\nabla''(u_\ell-u_\infty)\|^p_{L^p(\Omega_\ell;{\mathbb R}^{n-r})}.
$$
However, when $k>0$, we make crucial use of the uniform strict convexity to derive the estimate.

Let us close this section with an estimate similar to that obtained in Lemma~\ref{lem2}. Recall that $u_\ell$ is a minimizer on $\Omega_\ell$, whereas the following estimate is on $\Omega_{\ell_0}$. {See \cite{Chipot-Savitska} for a very similar argument.}

\begin{lemma}\label{estimation sur lzero}
There exist constants $\bar\ell$ and $c_4$, independent of $\ell$, such that for all $\bar \ell\le\ell_0\le \ell$,
\begin{equation}\label{estmtul zero}
\int_{\Omega_{\ell_0}}|\nabla u_{\ell}|^p\,dx\le c_4 \ell_0^r.
\end{equation}
\end{lemma}

\noindent{\bf Proof.}\enspace{} Let $1\le t\le \ell-1$ and set $\rho_t=\rho_{t+1,t}$. We take $v_{t,\ell}=(1-\rho_t)u_{\ell}$ as a test-function in problem~\eqref{pl}. This test-function is equal to $u_\ell$ ``far away'' and is $0$ in $\Omega_t$. We obtain
\begin{align*}
\int_{\Omega_{\ell}} F(\nabla u_{\ell})\,dx&\le \int_{\Omega_{\ell}}\bigl( F(\nabla v_{t,\ell})-f''(v_{t,\ell}-u_\ell)\bigr)\,dx\\
&=\int_{\Omega_t}F(0)\,dx+\int_{\Omega_t}f''u_\ell\,dx+
\int_{\Omega_{t+1}\setminus\Omega_t}\bigl[ F(\nabla v_{t,\ell})+f''\rho_tu_\ell)\bigr]\,dx\\
&\qquad\qquad\qquad
+\int_{\Omega_{\ell}\setminus\Omega_{t+1}} F(\nabla u_{\ell})\,dx.
\end{align*}
Therefore, we see that
$$
\int_{\Omega_{t+1}} F(\nabla u_{\ell})\,dx\le At^r+\int_{\Omega_{t+1}}\nu_tf''u_\ell\,dx+\int_{\Omega_{t+1}\setminus\Omega_t} F(\nabla v_{t,\ell})\,dx,
$$
with $A=F(0){\cal L}^r(\omega'){\cal L}^{n-r}(\omega'')$ and $\nu_t=\mathbf{1}_{\Omega_t}+\rho_t\mathbf{1}_{\Omega_{t+1}\setminus\Omega_t}$.

By the coerciveness and growth hypotheses \eqref{growth1}, we infer that
$$
\lambda\int_{\Omega_{t+1}} |\nabla u_{\ell}|^p\,dx\le Bt^r+\int_{\Omega_{t+1}}|f''u_\ell|\,dx+\Lambda\int_{\Omega_{t+1}\setminus\Omega_t} |\nabla v_{t,\ell}|^p\,dx,
$$
for some constant $B$, since $0\le\nu_t\le 1$ and the Lebesgue measure of $\Omega_{t+1}\setminus\Omega_t$ is of the order of $t^{r-1}$.

In $\Omega_{t+1}\setminus\Omega_t$, we have
$$|\nabla v_{t,\ell}|^p=|(1-\rho_t)\nabla u_{\ell}-u_\ell\nabla\rho_t|^p\le 2^{p-1}\bigl(|\nabla u_\ell|^p+K^p|u_\ell|^p\bigr).$$
Clearly, estimate \eqref{poincare} is also valid on $\Omega_{t+1}\setminus\Omega_t$, thus,
$$
\int_{\Omega_{t+1}\setminus\Omega_t}|u_\ell|^p\,dx
\le c_1^p\int_{\Omega_{t+1}\setminus\Omega_t}|\nabla''u_\ell|^p\,dx\le c_1^p\int_{\Omega_{t+1}\setminus\Omega_t}|\nabla u_\ell|^p\,dx,
$$
so that
$$\int_{\Omega_{t+1}\setminus\Omega_t} |\nabla v_{t,\ell}|^p\,dx\le 2^{p-1}(1+c_1^pK^p)\int_{\Omega_{t+1}\setminus\Omega_t} |\nabla u_{\ell}|^p\,dx.$$
Furthermore,
\begin{align*}
\int_{\Omega_{t+1}}|f''u_\ell|\,dx&\le\frac{\varepsilon}{p}\int_{\Omega_{t+1}}|u_\ell|^p\,dx +\frac{(t+1)^r}{\varepsilon^{p'/p}p'}{\cal L}^r(\omega')\|f''\|_{L^{p'}(\omega'')}^{p'}\\
&\le\frac{\varepsilon c_1^p}{p}\int_{\Omega_{t+1}}|\nabla u_\ell|^p\,dx +\frac{C}{\varepsilon^{p'/p}}t^r,
\end{align*}
with $\varepsilon>0$ to be chosen afterwards.

Let us set
$$h(t)=\int_{\Omega_t}|\nabla u_\ell|^p\,dx.$$
Putting all the above estimates together, it follows that
\begin{equation}
\lambda' h(t+1)\le E\bigl(h(t+1)-h(t)\bigr)+Dt^r,\label{une estimation de plus}
\end{equation}
with $\lambda'=\lambda-\frac{\varepsilon c_1^p}{p}$, $D=B+\frac{C}{\varepsilon^{p'/p}}$ and $E=2^{p-1}\Lambda(1+c_1^pK^p)$. We now pick $\varepsilon$ in such a way that $\lambda'>0$. Inequality \eqref{une estimation de plus} may be rewritten as
\begin{equation}\label{encore une recurrence a venir}
h(t)\le \theta h(t+1)+Ht^r,
\end{equation}
where $\theta=1-\frac{\lambda'}{E}\in {]}0,1[$ and $H=\frac D{E}$ depend neither on $t$ nor on $\ell$. Iterating inequality \eqref{encore une recurrence a venir}, we see that for $n=\lfloor\ell-t\rfloor$, we have
\begin{equation}\label{apres recurrence}
h(t)\le \theta^nh(t+n)+H\sum_{m=0}^{n-1}(t+m)^r\theta^m.
\end{equation}

Let us now set $t=\ell_0$. We have $h(\ell_0+\lfloor\ell-\ell_0\rfloor)\le h(\ell)\le c_2\ell^r$ by Lemma \ref{lem2}. Hence
$$\theta^{\lfloor\ell-\ell_0\rfloor}h(t+\lfloor\ell-\ell_0\rfloor)\le c_2 \theta^{\lfloor\ell-\ell_0\rfloor}\ell^r\le c_2 \theta^{\ell-\ell_0-1}\ell^r.$$
Now, for $\ell_0\ge-\frac{r}{\ln\theta}$, the function in the right-hand side is decreasing, hence maximum for $\ell=\ell_0$. Therefore,
$$\theta^{\lfloor\ell-\ell_0\rfloor}h(t+\lfloor\ell-\ell_0\rfloor)\le \frac{c_2}{\theta} \ell_0^r,$$
for $\ell\ge\ell_0\ge-\frac{r}{\ln\theta}$. Moreover, for $\ell_0\ge 1$,
\begin{multline*}
\sum_{m=0}^{\lfloor\ell-\ell_0\rfloor-1}(\ell_0+m)^r\theta^m=\ell_0^r\sum_{m=0}^{\lfloor\ell-\ell_0\rfloor-1}\Bigl(1+\frac m{\ell_0}\Bigr)^r\theta^m\\
\le \ell_0^r\sum_{m=0}^{\lfloor\ell-\ell_0\rfloor-1}(1+m)^r\theta^m\le \frac{\sum_{m=1}^{+\infty}m^r\theta^m}{\theta}\ell_0^r,
\end{multline*}
which completes the proof with  $\bar \ell=\max\bigl(1,-\frac{r}{\ln\theta}\bigr)$.\hfill$\square$\par\medbreak

We now turn to the convergence results. As a consequence of Lemma \ref{estimation sur lzero}, we have, without any restriction on $r$ with respect to $p$ and $k$,

\begin{theorem}\label{sous suite}
There exists a subsequence $\ell\to+\infty$ and a function $u^*\in W^{1,p}_{\rm loc}(\Omega_\infty)$ such that, for all $\ell_0$,
\begin{equation}\label{conv faible}
u_{\ell|\Omega_{\ell_0}}\rightharpoonup u^*_{|\Omega_{\ell_0}}\text{ weakly in }W^{1,p}(\Omega_{\ell_0}).
\end{equation}
Moreover, $u^*=0$ on $\partial\Omega_\infty$.
\end{theorem}

Note that the weak convergence above implies that $u_{\ell}\rightharpoonup u^*$ weakly in $W^{1,p}_{\rm loc}(\Omega_\infty)$. We will sometimes omit the restriction notation in the sequel when unnecessary.

\medskip
\noindent{\bf Proof.}\enspace{} By estimates \eqref{poincare} and \eqref{estmtul zero}, for all $n\in{\mathbb N}^*$, $u_\ell$ is bounded in $W^{1,p}(\Omega_n)$. Using the diagonal procedure, we thus construct a sequence $\ell_n$ such that for all $m$, $u_{\ell_n|\Omega_m}\rightharpoonup u^*_m$  weakly in $W^{1,p}(\Omega_m)$, with $u_m=0$ on $\omega'_m\times\partial\omega''$. Now, since $\Omega_m\subset\Omega_{m'}$ as soon as $m\le m'$, it follows that $u^*_m=u^*_{m'|\Omega_m}$, so that we have constructed a single limit function $u^*$ in the desired class. Furthermore, for all $\ell_0$, if we choose an integer $m\ge \ell_0$, we see that convergence~\eqref{conv faible} holds true.\hfill$\square$\par\medbreak

In the sequel, we will always consider a weakly convergent subsequence $u_\ell$ in the sense of Theorem \ref{sous suite}.

\section{Identification of the limit when $\ell\to+\infty$}
In this section, we do not make any further use of assumption \eqref{uniformite stricte} of uniform strict convexity of $F''$, other than the fact that we used it to establish Theorem \ref{estimation de base}.\footnote{Keep in mind that this hypothesis is void for $k=0$ anyway.} The results will only hold for values of $r$ small enough depending on $p$. We let $\Omega_\infty={\mathbb R}^r\times\omega''$.

Let us first show that the asymptotic behavior of $u_\ell$ is independent of the elongated dimension if $r$ is small enough.

\begin{theorem}\label{limite allongee}
Assume that $r<p$ if $k=0$, or that $r<kp/(p-k)$ if $0<k<p$. Then
we have $\nabla'u^*=0$ and $u^*$ may be identified with a function in the $x''$ variable only, still denoted $u^*$, which belongs to $W^{1,p}_0(\omega'')$.
\end{theorem}

\noindent{\bf Proof.}\enspace{} By estimates \eqref{estmtul}  and  \eqref{Caccio usc} and the triangle inequality, it follows that
\begin{equation}\label{gradient' tend vers zero}
\|\nabla'u_\ell\|^p_{L^p(\Omega_{\ell_0};{\mathbb R}^r)}\le C\biggl(\frac{\delta}{(\ell-\ell_0)^p}+\frac{k}{(\ell-\ell_0)^{kp/(p-k)}}\biggr)\ell^{r}\to 0
\end{equation}
when $\ell\to+\infty$ with $\ell_0$ fixed. Indeed, when $0<k\le p/2$, we actually have $\frac{kp}{p-k}\le p$ and since $\ell\to+\infty$, the first term in the right hand side of estimate \eqref{gradient' tend vers zero} is bounded from above by the second term. 

Now $\nabla'u_\ell\rightharpoonup\nabla'u^*$ weakly in $L^{p}_{\rm loc}(\Omega_\infty)$, hence we see that $\nabla'u^*=0$, which concludes the proof of the Theorem.\hfill$\square$\par\medbreak

In order to get a feeling of what Theorem \ref{limite allongee} says, let us look at a few examples. For the Laplacian, we have $p=2$ and we can take $k=0$, which restricts this result to $r=1$ (see Section 5 for a more general result  with additional hypotheses, that applies in this case). For the $p$-Laplacian, $p>2$, we can take $k=2$ and the result is restricted to $r<2p/(p-2)$.    {This restriction for the $p$-Laplacian can already be found in \cite{Xie}}. Note that $r=1$ and $r=2$ are allowed for any value of $p$. This is not optimal in this particular case, since it is known that $\ell\to+\infty$ convergence holds without restriction on the dimension with respect to $p$, see \cite{Chipot-Xie}.

Let us now identify the limit function. We first need another estimate.

\begin{lemma}\label{estimation sur une tranche}
There exists a constant $c_5$ such that for all $t\le s$,
\begin{equation}\label{estime tranche}
\limsup_{\ell\to+\infty}\int_{\Omega_s\setminus\Omega_t}|\nabla u_{\ell}|^p\,dx\le c_5 (s^r-t^r).
\end{equation}
\end{lemma}

\noindent{\bf Proof.}\enspace{} We may assume that $t>0$, since the case $t=0$ is already covered by Lemma \ref{estimation sur lzero}. We use here De Giorgi's classical slicing trick. Let $n$ be an integer large enough so that $0\le t-\frac1n<s+\frac1n\le\ell$. For each integer $m$,  $1\le m\le n$, we consider the cut-off function
$$\chi_{m,n}(x')=\rho_{s+\frac{m}{n^2},s+\frac{m-1}{n^2}}(x')\Bigl(1-\rho_{t-\frac{m-1}{n^2},t-\frac{m}{n^2}}(x')\Bigr).$$
This cut-off function takes its values in $[0,1]$, it is $0$ whenever $g(x')\ge s+\frac{m}{n^2}$ or $g(x')\le t-\frac{m}{n^2}$, it is $1$ for $t-\frac{m-1}{n^2}\le g(x')\le s+\frac{m-1}{n^2}$, and $|\nabla\chi_{m,n}|\le Kn^2$. Let us call $S_{m,n}$ the slice where $0<\chi_{m,n}(x')<1$. We observe that
$$\bigcup_{m=1}^n \overline{S_{m,n}}=\overline{\Omega_{s+\frac1n}\setminus\Omega_s}\bigcup \overline{\Omega_{t}\setminus\Omega_{t-\frac1n}}
\subset \Omega_{s+1},$$
and that $S_{m,n}\cap S_{m',n}=\emptyset$ when $m\neq m'$.

Let us consider the test-function $v_{\ell,m,n}=(1-\chi_{m,n})u_\ell+\chi_{m,n}u^*$. The minimization problem yields the estimate
\begin{align*}
\int_{\Omega_\ell}F(\nabla u_\ell)\,dx&\le\int_{\Omega_\ell}F(\nabla v_{\ell,m,n})\,dx-\int_{\Omega_\ell}f''\chi_{m,n}(u^*-u_\ell)\,dx\\
&=\int_{\Omega_\ell}F(\nabla v_{\ell,m,n})\,dx-\int_{\Omega_{s+1}}f''\chi_{m,n}(u^*-u_\ell)\,dx.
\end{align*}
Taking into account the specific form of the cut-off function, this implies that
\begin{align}
\nonumber\int_{\Omega_s\setminus\Omega_t}F(\nabla u_\ell)\,dx&\le \int_{\Omega_{s+\frac{m}{n^2}}\setminus\Omega_{t-\frac{m}{n^2}}}F(\nabla u_\ell)\,dx\\
\nonumber&\le\int_{\Omega_{s+\frac{m}{n^2}}\setminus\Omega_{t-\frac{m}{n^2}}}F(\nabla v_{\ell,m,n})\,dx-\int_{\Omega_{s+1}}f''\chi_{m,n}(u^*-u_\ell)\,dx\\
\nonumber&\le\int_{S_{m,n}}F(\nabla v_{\ell,m,n})\,dx+\int_{\Omega_{s+\frac{m-1}{n^2}}\setminus\Omega_{t-\frac{m-1}{n^2}}}F(\nabla u^*)\,dx\\
&\qquad\qquad-\int_{\Omega_{s+1}}f''\chi_{m,n}(u^*-u_\ell)\,dx.\label{slicing1}
\end{align}
Let us estimate each term in the right-hand side separately. First of all, we have
\begin{equation}\label{slicing2}
\Bigl|\int_{\Omega_{s+1}}f''\chi_{m,n}(u^*-u_\ell)\,dx\Bigr|\le A^{1/p'}(s+1)^{r/p'}\|f''\|_{L^{p'}(\omega'')}\|u^*-u_\ell\|_{L^p(\Omega_{s+1})},
\end{equation}
with $A=\mathcal{L}^r(\omega')$.
Secondly, we see that
\begin{equation}\label{slicing3}
\biggl|\int_{\Omega_{s+\frac{m-1}{n^2}}\setminus\Omega_{t-\frac{m-1}{n^2}}}F(\nabla u^*)\,dx\biggr|\le A\Bigl(\Bigl(s+\frac1n\Bigr)^r-
\Bigl(t-\frac1n\Bigr)^r\Bigr)\|F''(\nabla'' u^*)\|_{L^1(\omega'')}.
\end{equation}

We now come to the slicing argument stricto sensu. By the growth estimate \eqref{growth1}, we have
\begin{multline}\label{slicing4}
\int_{S_{m,n}}F(\nabla v_{\ell,m,n})\,dx\le2^{p-1}\Lambda\Bigl(\int_{S_{m,n}}\bigl(|\nabla u_\ell|^p+|\nabla u^*|^p+1\bigr)\,dx\\+K^pn^{2p}
\int_{S_{m,n}}|u^*- u_\ell|^p\,dx\Bigr).
\end{multline}
The only term that causes a difficulty is the last term coming from $\nabla\chi_{m,n}$.
We now plug estimates \eqref{slicing2}, \eqref{slicing3} and \eqref{slicing4} into the right-hand side of estimate \eqref{slicing1}, sum for $m=1$ to $n$ and divide the result by $n$. Observing that the sum of integrals over the slices $S_{m,n}$ gives rise to integrals over the union of all slices, which is included in $\Omega_{s+1}$, this yields
\begin{align}
\int_{\Omega_s\setminus\Omega_t}F(\nabla u_\ell)\,dx&\le A^{1/p'}(s+1)^{r/p'}\|f''\|_{L^{p'}(\omega'')}\|u^*-u_\ell\|_{L^p(\Omega_{s+1})}\\
&\qquad +A\Bigl(\Bigl(s+\frac1n\Bigr)^r-
\Bigl(t-\frac1n\Bigr)^r\Bigr)\|F''(\nabla'' u^*)\|_{L^1(\omega'')}\\
&+\frac{2^p\Lambda c_4}n(s+1)^r+2^{p-1}\Lambda K^p n^{2p-1}\|u^*-u_\ell\|^p_{L^p(\Omega_{s+1})}.
\end{align}

We first let $\ell\to+\infty$. Due to the Rellich-Kondra\v sov theorem, $\|u^*-u_\ell\|_{L^p(\Omega_{s+1})}\to 0$ and it follows from the coerciveness  estimate that
\begin{multline*}
\limsup_{\ell\to+\infty}\int_{\Omega_s\setminus\Omega_t}|\nabla u_{\ell}|^p\,dx\le \frac A\lambda\Bigl(\Bigl(s+\frac1n\Bigr)^r-
\Bigl(t-\frac1n\Bigr)^r\Bigr)\|F''(\nabla'' u^*)\|_{L^1(\omega'')}\\+\frac{2^p\Lambda c_4}{n\lambda}(s+1)^r.
\end{multline*}
We finally let $n\to+\infty$ to obtain the result with $c_5=\frac A\lambda\|F''(\nabla'' u^*)\|_{L^1(\omega'')}$.\hfill$\square$\par\medbreak

We now are in a position to prove the main result of this section.

\begin{theorem}\label{resultat principal faible}
The function $u^*$ is a minimizer of problem ${\cal P}_\infty$.
\end{theorem}

\noindent{\bf Proof.}\enspace{}   Let $z\in W^{1,p}_0(\omega'')$ be arbitrary. We use the test function
$v_\ell=(1-\rho_{t})u_\ell+\rho_{t}z$, with  $\rho_t=\rho_{t+1,t}$, so that $v_\ell=u_\ell$ on $\Omega_\ell\setminus\Omega_{t+1}$ and $v_\ell=z$ on $\Omega_t$. We thus have
\begin{equation}\label{on va finir par y arriver}
\int_{\Omega_{t+1}}[F(\nabla u_\ell)-f''u_\ell]\,dx\le \int_{\Omega_{t+1}\setminus\Omega_t}[F(\nabla v_\ell)-f''v_\ell]\,dx
+\int_{\Omega_{t}}[F(\nabla z)-f''z]\,dx.
\end{equation}

It follows from Lemma \ref{estimation sur une tranche} that 
$$\limsup_{\ell\to+\infty}\Bigl|\int_{\Omega_{t+1}\setminus\Omega_t}[F(\nabla v_\ell)-f''v_\ell]\,dx\Bigr|\le C(t+1)^{r-1}$$
for some constant $C$ independent of $\ell$ and $t$. The left-hand side of estimate \eqref{on va finir par y arriver} is weakly lower-semicontinuous, hence, letting $\ell\to+\infty$, we obtain
\begin{align*}
(t+1)^r\mathcal{L}^r(\omega') \int_{\omega''}[F(\nabla u^*)-f''u^*]\,dx'&\le C(t+1)^{r-1}\\
&\qquad+t^r\mathcal{L}^r(\omega') \int_{\omega''}[F(\nabla z)-f''z]\,dx'
\end{align*}
and the result follows from letting $t\to+\infty$, since $F(\nabla u^*)=F''(\nabla'' u^*)$ and $F(\nabla z)=F''(\nabla'' z)$.\hfill$\square$\par\medbreak

We now apply a classical trick to obtain strong convergence when $F''$ is strictly convex. Of course, when $k>0$, this is already the case by assumption \eqref{uniformite stricte}. Strict convexity is only a new assumption if $k=0$. In this case, the solution $u_\infty$ of the limit problem is unique and this uniqueness implies the weak convergence of the whole family $u_\ell$.

\begin{theorem}\label{resultat principal fort}
Assume that $F''$ is strictly convex. Then $u^*=u_\infty$ and 
$u_{\ell}\to u_\infty$ strongly in $W^{1,p}(\Omega_{\ell_0})$ for all $\ell_0$.
\end{theorem}

We recall the following two lemmas that can be found \emph{e.g.} in \cite{Ball-Marsden}.

\begin{lemma}\label{BM1} Let $F\colon {\mathbb R}^M\to {\mathbb R}$ be strictly convex. Let $\mu\in{]}0,1[$ and $a_j,a\in {\mathbb R}^M$ such that
$$
\mu F(a_j)+(1-\mu)F(a)-F(\mu a_j+(1-\mu) a)\to 0\text{ as }j\to+\infty.
$$
Then $a_j\to a$.
\end{lemma}

The second lemma is a slight variation on Fatou's lemma.
\begin{lemma}\label{BM2} Let $F_j,F,H_j,H\in L^1(\Omega)$ with $F_j\ge H_j\ge 0$ for all $j$, $F_j\to F$ and $H_j\to H$ a.e., and $\int_\Omega F_j\,dx\to\int_\Omega F\,dx$. Then
$$
\int_\Omega H_j\,dx\to\int_\Omega H\,dx.
$$
\end{lemma}

\noindent{\bf Proof of Theorem \ref{resultat principal fort}.} We already know that $\nabla' u_\ell\to 0=\nabla' u^*$ strongly in $L^p(\Omega_t)$ by estimate \eqref{gradient' tend vers zero}. We thus just have to prove the strong convergence of $\nabla'' u_\ell$.

We use a similar slicing as before, with the test-functions $\rho_{t+\frac m{n^2},t+\frac {m-1}{n^2}}u_\ell$ for $n$ large enough, $1\le m\le n$. Skipping the details, this slicing implies that
$$\limsup_{\ell\to+\infty}\int_{\Omega_t}F(\nabla u_\ell)\,dx\le \int_{\Omega_t}F(\nabla u^*)\,dx.$$

On the other hand, for almost all $x'$, the function $u_{x',\ell}\colon x''\mapsto u_\ell(x',x'')$ is an admissible test-function for the limit problem, so that
$$\int_{\omega''}[F''(\nabla''u^*)-f''u^*]\,dx''\le \int_{\omega''}[F''(\nabla''u_{x',\ell})-f''u_{x',\ell}]\,dx''.$$
We integrate this inequality with respect to $x'\in t\omega'$ and obtain
$$\int_{\Omega_t}[F''(\nabla''u^*)-f''u^*]\,dx\le \int_{\Omega_t}[F''(\nabla''u_\ell)-f''u_\ell]\,dx.$$
We now let $\ell\to+\infty$, which yields
$$\int_{\Omega_t}F''(\nabla''u^*)\,dx\le \liminf_{\ell\to+\infty}\int_{\Omega_t}F''(\nabla''u_\ell)\,dx.$$
By hypothesis \eqref{growth3}, $G\ge 0$, which implies that $F''(\xi'')\le F(\xi',\xi'')$ for any $\xi'$. It follows that
\begin{equation}\label{func-conv}
  \int_{\Omega_t}F''(\nabla''u_\ell)\,dx\to \int_{\Omega_t}F''(\nabla''u^*)\,dx
\end{equation}
when $\ell\to+\infty$, since $F''(\nabla''u^*)=F(\nabla u^*)$.

Let us pick $\mu\in{]}0,1[$ and set
$$g_\ell=\mu F''(\nabla'' u_\ell)+(1-\mu)F''(\nabla'' u^*)-F''(\mu \nabla'' u_\ell+(1-\mu)\nabla'' u^*).
$$
By weak lower semicontinuity, it is clear that
$$\liminf_{\ell\to+\infty}\int_{\Omega_t}F''(\mu \nabla'' u_\ell+(1-\mu)\nabla'' u^*)\,dx
\ge\int_{\Omega_t}F''(\nabla'' u^*)\,dx.$$
Therefore
$$0\le \limsup_{\ell\to+\infty}\int_{\Omega_t}g_\ell\,dx\le \int_{\Omega_t}F''(\nabla'' u^*)\,dx-\int_{\Omega_t}F''(\nabla'' u^*)\,dx=0,$$
so that $g_\ell\to 0$ a.e.\ (up to a subsequence).
We then apply Lemma \ref{BM1} to deduce that $\nabla''u_\ell\to\nabla''u^*$ a.e.\ up to that same subsequence.

We now let
$$H_\ell=|\nabla''u_\ell-\nabla''u^*|^p\le 2^{p-1}(F''(\nabla'' u_\ell)+|\nabla''u^*|^p)=F_\ell,$$
and invoke Lemma \ref{BM2} and \eqref{func-conv} to obtain the result for $\ell_0=t$. To conclude for all $\ell_0$, we use the diagonal process.\hfill$\square$\par\medbreak

\section{Convergence rates}
In the previous section, we obtained convergence results without taking advantage of the term involving $k$ in the left-hand side of estimate \eqref{Caccio usc}. This makes them valid in particular for $k=0$ without strict or uniform strict convexity. It should however be clear that for $k>0$, the term in question can be used to obtain a much shorter convergence proof with convergence rate, which we do not detail here. More precisely,

\begin{theorem}\label{limite allongee usc}
Under the previous hypotheses with $0<k<p$ and $r<kp/(p-k)$, we have
$$ \|u_{\ell}-u_{\infty}\|^p_{W^{1,p}(\Omega_{\ell_0})}\leq C\ell^{r-\frac{kp}{p-k}}.$$
\end{theorem}

The proof is a direct consequence of Corollary \ref{Caccio global} and Lemma \ref{lem2}.

In any case, the estimates do not seem to allow a convergence proof without any restriction on $r$ with respect to $p$ in all generality, whereas it is known in some cases, for instance in the case of the Laplacian, that convergence holds true for all values of $r$.

In order to partially overcome these shortcomings, we assume now that $k=0$ and that $F''$ is uniformly strictly convex in the sense that \begin{equation}\label{uniformite stricte k=0}
F''(\theta\xi''+\mu\zeta'')\le \theta F''(\xi'')+\mu F''(\zeta'')-\beta\theta\mu(\theta^{p-1}+\mu^{p-1})|\xi''-\zeta''|^p,
\end{equation}
for some $\beta>0$. Note that this is equivalent to allowing $k=p$ in hypotheses \eqref{growth3} and \eqref{uniformite stricte}. In some sense, $\frac{kp}{p-k}$ is then infinite and it is to be expected that there should be no restriction on the allowed dimensions $r$, plus faster than polynomial convergence. This is what we now proceed to show.

Under assumption \eqref{uniformite stricte k=0}, it is fairly clear that we still have an estimate similar to that of  Theorem \ref{estimation de base}, namely,
\begin{equation}\label{Caccio ter}
\|\nabla'u_\ell\|^p_{L^p(\Omega_{t})}+\|\nabla''(u_\ell-u_\infty)\|^p_{L^p(\Omega_t)}\le \frac{C}{(s-t)^p}\|\nabla''(u_\ell-u_\infty)\|^p_{L^p(\Omega_s\setminus\Omega_t)}.
\end{equation}

Let us thus prove that not only does convergence hold without restrictions on the elongated dimension $r$, but that it also occurs at an exponential rate. The extra control makes things actually much easier.

\begin{theorem}\label{th4}Under hypotheses \eqref{growth1}-\eqref{growth3} with $k=0$ and \eqref{uniformite stricte k=0}, then for all $r< n$ and all $\ell_0$, there exist constants $C$ and $\alpha>0$ independent of $\ell$ such that we have
$$||\nabla(u_{\ell}- u_{\infty})||_{L^p(\Omega_{\ell_0})}\leq C e^{-\alpha\ell}.$$
\end{theorem}

\noindent{\bf Proof.}\enspace{}
We take $s=t+1$ in estimate \eqref{Caccio ter}, which yields
\begin{align*}
\|\nabla'u_\ell\|^p_{L^p(\Omega_{t})}+\|\nabla''(u_\ell-u_\infty)\|^p_{L^p(\Omega_t)}&\le C\|\nabla''(u_\ell-u_\infty)\|^p_{L^p(\Omega_{t+1}\setminus\Omega_{t})}\\ &\leq
C\|\nabla'u_\ell\|^p_{L^p(\Omega_{t+1}\setminus\Omega_{t})}\\
&\qquad\qquad{}+
C\|\nabla''(u_\ell-u_\infty)\|^p_{L^p(\Omega_{t+1}\setminus\Omega_{t})}.
\end{align*}
Setting
 $$g(t)=\|\nabla'u_\ell\|_{L^p(\Omega_t)}^p+\|\nabla''(u_\ell-u_\infty)\|_{L^p(\Omega_t)}^p,$$
we have just shown that
$$
g(t)\le C(g(t+1)-g(t)),
$$
or in other words
\begin{equation}\label{la, ca change encore-new5}
g(t)\le \theta g(t+1)
\end{equation}
with $\theta=\frac{C}{1+C}\in{}]0,1[$.

We iterate inequality \eqref{la, ca change encore-new5} using the sequence $t_n=n+\ell_0$, $n=0,\ldots,\lfloor \ell-\ell_0\rfloor$. Obviously
$$g(\ell_0)=g(t_0)\leq \theta^ng(t_n)$$
for all such $n$, and in particular for the last one,
$$g(\ell_0)\leq\theta^{\lfloor \ell-\ell_0\rfloor}g(t_{\lfloor \ell-\ell_0\rfloor})\le 
\theta^{\ell-\ell_0-1}g(\ell)\leq C\theta^{-\ell_0-1}e^{\ell\ln \theta}\ell^r,$$ 
with $\ln\theta<0$. Now, for all $r$, we can pick $\alpha>0$ such that $\ln\theta<-p\alpha<0$ and $e^{\ell\ln \theta}\ell^r\le e^{-p\alpha \ell}$ for $\ell$ large enough, which completes the proof since $\nabla' u_\infty=0$.\hfill$\square$\par\medbreak

Theorem \ref{th4} applies to energies of the form $F(\xi)=F'(\xi')+F''(\xi'')$, for instance. We recover in particular the known result for the case of the 2-Laplacian. {See also the monograph \cite{[C2new]} for exponential estimates in this context.}

\section{Extension to the vectorial case}

We have written everything so far in the context of a scalar problem, \emph{i.e.}, the functions $u_\ell$ are scalar-valued. All previous developments  only made use of the minimization problem, under various convexity assumptions. Now clearly, absolutely nothing is changed if we consider instead vector-valued problems in the calculus of variations, with functions $u_\ell$ taking their values in some ${\mathbb R}^N$, if the energies are supposed to satisfy the same growth, coercivity and convexity assumptions as before, and the same convergence results hold true.

Unfortunately, in the vectorial case of the calculus of variations, the relevant condition that guarantees lower-semicontinuity of the energy functional is not convexity, but much weaker conditions such as quasiconvexity, or in the case of energies that can take the value $+\infty$, as is the case in nonlinear elasticity, polyconvexity, see \cite{[D]}. Indeed, convexity is not suitable in nonlinear elasticity for well-known modeling reasons. This explains why we have striven to use as little convexity as possible (in some sense) at any given point in the sequence of arguments. This comment should however be mitigated by the fact that some instances of our uses of convexity will also work with rank-1-convexity, which is a reasonable assumption in the vectorial case. There are also notions of strict uniform quasiconvexity that may apply, see \cite{Evans}.

The fact that the Euler-Lagrange equation is not available in nonlinear elasticity is also an incentive to try and only use the minimization problem.
Now, it is at this point unclear to us how to attack the elongation problem in such nonconvex vectorial cases, since we still heavily rely on (strict uniform) convexity at crucial points of the proofs. Moreover, the Dirichlet boundary condition considered here is not necessarily the most interesting one in the context of nonlinear elasticity,  in particular if we have the Saint Venant principle in mind.

Even the potential limit problem is not so clear. In another dimension reduction context, when considering a body whose thickness goes to zero, and with different boundary conditions, it can be seen that quasiconvexity is not conserved through an ``algebraic'' formula of the kind found here, and that a relaxation step is necessary, see for instance \cite{[L-R]}. Physically, this due to the possibility of crumpling such a thin body. A similar phenomenon may quite possibly happen here, but maybe not in the same fashion.

To the best of our knowledge, the nonconvex vectorial case remains open.

\end{document}